\def\intav{-\hskip -1.1em\int}
\input amstex
\magnification=\magstep1
\documentstyle{amsppt}
\NoRunningHeads
\pageheight{43pc}
\pagewidth{28pc}

\topmatter
\title Boundary limits for bounded quasiregular mappings \\
\endtitle
\author Bao Qin Li$^\ast$\\
	Department of Mathematics\\
	Florida International University\\
	Miami, FL 33199, \\
        \\and \\
	\\
	Enrique Villamor\\
	Department of Mathematics\\
	Florida International University\\
	Miami, FL 33199\\
\endauthor
\thanks
$^\ast$ Research partially supported by NSF grant DMS-9101864
\endthanks



\abstract
	In this paper we establish results on the existence
	of nontangential limits for weighted $\Cal A$-harmonic functions
	in the weighted Sobolev space $W_w^{1,q}(\Bbb B^n)$, for some $q>1$ and 
	$w$ in the Muckenhoupt $A_q$ class,
	where $\Bbb B^n$ is the unit ball in $\Bbb R^n$. These results generalize the ones in section \S3 
	of [KMV], where the weight was identically equal to one.
	Weighted $\Cal A$-harmonic functions are weak solutions of the partial differential 
	equation
	$$\text{div}(\Cal A(x,\nabla u))=0,$$
	where $\alpha\,w(x)\,|\xi|^{q} \le \langle \Cal A(x,\xi),\xi \rangle\le 
	\beta\,w(x)\, |\xi|^{q}$ for some fixed $q\in (1,\infty)$, where $0<\alpha\leq \beta<\infty$, and 
	$w(x)$ is a $q$-admissible weight as in Chapter 1 in [HKM].\par
	Later, we apply these results to 
	improve on results of Koskela, Manfredi and Villamor [KMV] and Martio and Srebro [MS] 
	on the existence of radial 
	limits for bounded quasiregular mappings in the unit ball of $\Bbb R^n$ with some 
	growth restriction on their multiplicity function.   
\endabstract
\endtopmatter
\document
\beginsection{$\S$1. Introduction.}

In this paper we study weak solutions of the partial differential equation
$$\text{div}(\Cal A(x,\nabla u))=0,\tag 1.1$$ 
where $\Cal A \colon \Bbb R^n\times \Bbb R^n\to\Bbb R^n$ is a mapping satisfying the following assumptions 
for some constants $0<\alpha \leq \beta<\infty$:\par
$$\text{the mapping } x\to \Cal A(x, \xi)\text{ is measurable for all }\xi\in \Bbb R^n \text{ and}$$
$$\tag 1.2$$
$$\text{the mapping } \xi\to \Cal A(x, \xi)\text{ is continuous for a.e. } x\in \Bbb R^n;$$
for all $\xi\in \Bbb R^n$ and a.~e.~ $x\in \Bbb R^n$
$$\langle \Cal A(x,\xi),\xi\rangle\geq \alpha \,w(x)\, |\xi|^q\tag 1.3$$
$$|\Cal A(x,\xi)|\leq \beta\,w(x)\,|\xi|^q\tag 1.4$$
where $1<q<\infty$;
$$\langle(\Cal A(x,\xi_1)-\Cal A(x,\xi_2)), (\xi_1-\xi_2)\rangle >0 \tag 1.5$$
whenever $\xi_1,\,\xi_2\in \Bbb R^n$, $\xi_1\neq\xi_2$; and
$$\Cal A(x,\lambda\,\xi)=\lambda\, |\lambda|^{q-2}\, \Cal A(x,\xi)\tag 1.6$$
whenever $\lambda\in \Bbb R$, $\lambda\neq 0$.
More generally we could have replaced in (1.3) and (1.4) $\beta$ by a function $\beta(x)$ with the condition that is 
bounded and $\alpha$ by a function $\alpha(x)$ asking that $\alpha(x)>0$ a.~e.~ $x$. 
Instead, we will consider the uniformly elliptic case and general $1<q<\infty$. 
Here we assume that $w(x)$ is a $q$-admissible nonnegative weight as defined in Chapter 1 of [HKM].

Solutions of (1.1) are called weighted $\Cal A$-harmonic functions. The  
prototype of these equations is the weighted $p$-Laplace equation
$$\Delta_{q,w}u=\text{div}(w(x)\,|\nabla u|^{q-2}\nabla u)=0.$$ \par 

In this note we present generalizations of a number of theorems on 
the existence of nontangential limits of  
weak solutions of (1.1)  in a ball $\Bbb B$ with finite weighted $q$-Dirichlet integral
$$\int_{\Bbb B} |\nabla u|^q \,w(x)\, dx<\infty,$$
for some $q>1$ and $w$ a nonnegative weight in the Muckenhoupt $A_q$ class. 

Notice that for $q>n-1$ monotonicity of the functions in the weighted 
Sobolev class is all it is needed as is shown in [MV1]. It is well known that solutions of (1.1) 
satisfying conditions (1.2) through (1.6) are monotone. Hence, the real contribution of this note for solutions of (1.1) 
satisfying conditions (1.2) through (1.6) is in the range $1<q\leq n-1$.\par
The weighted Sobolev class $W^{1,q}(\Bbb B^n;w)$ is defined in [HKM, Chapter 1]. 
It consists of functions
$u\colon \Bbb B^n \to \Bbb R^n$ that have first distributional 
derivatives $\nabla u$
such that $$\int_{\Bbb B^n} \left(|u(x)|^q+|\nabla u(x) |^q\right
)\,\,w(x)\,\,dx<\infty.$$ 
The weighted $q$-capacity we will be using throughout this paper is the 
relative first order variational $(q,w)$-capacity [HKM, Chapter 2]. Let us recall for the sake of 
completeness the definitions of monotone functions and of Muckenhoupt $A_q$ weights.\par

\proclaim{Definition 1.7}
	Let $\Omega\subset \Bbb R^n$ be an open set. A continuous function 
$u\colon \Omega\to \Bbb R$ is monotone, 
in the sense of Lebesgue, if 
$$\max_{\overline D} u(x)=\max_{\partial D} u(x)$$
and 
$$\min_{\overline D}u(x)=\min_{\partial D} u(x)$$ 
hold whenever $D$ is a domain with compact closure ${\overline 
D}\subset 
\Omega$.
\endproclaim

 \proclaim{Definition 1.8}
Let $q>1$ and $w\in L_{\text{loc}}^1(\Bbb R^n)$. We say that $w\in 
A_q$, 
if there exists a constant $C$ such that
$$\sup_{ B} \,\,\biggl(\intav_{ B} w(y)\,dy\biggr)\,\,\,
\biggl(\intav_{ B} w(y)^{{1\over{1-q}}}\, dy\biggr)^{q-1}<C$$
where the supremum is taken over all balls ${ B}\subset\Bbb R^n$.   
\endproclaim

Let us observe that if the weight 
$w$ is in $A_q$ it follows that it is $q$-admissible with the same index $q$, see Chapter 15 in [HKM]. \par
Our results extend the ones in section \S3 of [KMV], where the weight function $w$ was 
identically equal to one. We refer to the introduction of [KMV], for a 
historical chronology, background and references for these type of results. \par

In section \S2 we generalize the results in section \S3 of [KMV] for the weighted case. \par
In section \S3 we apply the results in section \S2 and in [MV1] to the components of 
bounded quasiregular mappings $f$ satisfying certain 
growth conditions on their multiplicity function $N(f,E)$ defined as follows: 
Let $E$ be a subset of the unit ball $\Bbb B^n$ in $\Bbb R^n$. We define  
$n(y;f,E)=\text{card} \lbrace x\in E\colon f(x)=y\rbrace$, and  
$N(f,E)=\sup_{y\in \Bbb R^n} n(y;f,E)$. $N(f,E)$ is called the {\it  
multiplicity  
function} of $f$. \par
Our main result appears in this section, and is the following:

\proclaim  {Theorem 3.11} Let $f$ be a bounded quasiregular mapping of $\Bbb B^n$, and suppose that, for some 
$0\leq a<n-1$,
$$N(f, B(0,r))\leq C\,(1-r)^{-a}$$ 
for all $0<r<1$. Then the set of points $x_0\in E\subset \partial \Bbb B^n(0,1)$ for which 
the nontangential limit of $f$ does not exist has Hausdorff dimension less 
than or equal to $a$, i.e. $\text{dim}_H(E)\leq a$.
\endproclaim

This theorem improves Theorem 4.1 in [KMV] where the bound on the Hausdorff dimension of $E$ was found to be 
${{n\,a}\over{1+a}}$. It is clear that for any $0<a<n-1$ we have that $a<{{n\,a}\over{1+a}}$. \par
In 1999, Martio and Srebro [MS] showed that for a bounded quasiregular mapping in the unit ball of 
$\Bbb R^n$ satisfying the same condition as in Theorem 3.11 for its multiplicity function , if $f$ is locally injective, 
then $\text{dim}_H(E)\leq a$, and that this estimate is best possible. For the sharpness of their 
estimate, they constructed for every $n\geq 3$ a sequence of numbers $s_m\in (0,n-1)$ and locally injective bounded 
quasiregular mappings $f_m$ in $\Bbb B^n$, $m=1,2, \ldots$ such that $f_m$ satisfies that 
$$N(f_m, B(0,r))\leq C\,(1-r)^{-s_m}$$
and $\text{dim}_H(E_m)=s_m$ with $\lim_{m\to\infty} s_m= n-1$, where $E_m$ is the set of points 
$x_0\in E_m\subset \partial \Bbb B^n(0,1)$ for which 
the nontangential limit of $f_m$ does not exist. These examples also show that our result is sharp. \par
Theorem 3.11 improves on Martio and Srebro's in that we do not assume in ours that the bounded 
quasiregular mappings are locally injective. As it will be shown at the end of section \S3, our result also holds 
with the same conclusion, for 
quasiregular mappings with the same restriction on the growth of their multiplicity function but not necessarily bounded, 
i.e. we can assume that $|f(x)|\leq C\, (1-|x|)^{-b}$ for some $0<b<\infty$ and Theroem 3.11 still holds with the same conclusion. 
Martio and Srebro's result only holds for locally injective and bounded 
qausiregular mappings with restricted growth in their multiplicity function. Thus, the natural question to ask, is 
whether or not Theorem 3.11 holds for $a=n-1$, that is:\par
Is it true that a bounded quasiregular mapping satisfying that
$$N(f, B(0,r))\leq C\,(1-r)^{-(n-1)}$$
has nontangential boundary limits everywhere on the 
boundary of the unit ball except possibly on a set of $(n-1)$ Hausdorff measure zero? \par

Recently Heinonen and Rickman [HR], have constructed examples in dimension $n=3$, 
of bounded quasiregular mappings, not necessarily 
locally injective, with no radial limits at points in sets of the boundary of the $\Bbb B^3$ with Hausdorff dimension 
arbitrarily close to 2.

\beginsection{ $\S$2.  Existence of Nontangential Limits  }

In this section we prove straightforward generalizations of results in section \S3 of [KMV] 
on boundary limits for weighted Dirichlet finite 
$\Cal  A$-harmonic functions. In [KMV] no weight was considered. Throughout this section we assume that
$\alpha(x)>\alpha>0$ for a.~e. $x$ and that our weights $w$ are in the $A_q$ 
Muckenhoupt class.\par
First we show that weighted Dirichlet finite $\Cal A$-harmonic functions $u,$ 
defined  
in the unit ball $\Bbb B^n$ of $\Bbb R^n$, have nontangential limits 
everywhere  
on the boundary of the unit ball except possibly on a set $E$ of weighted Bessel  
$B^w_{1,q}$-capacity zero, $1<q\leq n$. 
In this work we will use the weighted Bessel  
$B^w_{1,q}$-capacity for technical reasons. We refer the reader to the book by  
Ziemer [Z] for the definition and properties of the weighted Bessel capacity $B^w_{1,q}$, 
where the weight $w(x)\in A_q(\Bbb R^n)$.  
At this point, we would like to remark that all the $(q,w)$-capacities are  
equivalent in the sense that a set with one of the standard $(q,w)$-capacities  zero  
will have all the other $(q,w)$-capacities zero. Thus, in the rest of the paper we  
will say $(q,w)$-capacity zero without specifying the weighted capacity that we are using.  
The case $q>n$ is not interesting because then $u$ is continuous up to the  
boundary by the weighted Sobolev embedding theorem. Recall that a weighted $\Cal A$-harmonic  
function $u$ of $\Bbb B^n$ is continuous in $\Bbb B^n.$  \par

\proclaim{ Theorem 2.1}
Let $u$ be a weighted $\Cal A$-harmonic function in the unit ball $\Bbb B^n$ of $\Bbb  
R^n$ (no restriction on the type of $\Cal A$). If $\int_{\Bbb B^n} |\nabla  
u(x)|^q\,w(x)\,dx<\infty$ for some $1<q\le n$ and $w(x)\in A_q(\Bbb R^n)$, then the function $u$ has  
nontangential limits on all radii terminating outside a set of $(q,w)$-capacity  
zero.
\endproclaim
\noindent Theorem 2.1 extends Theorem 3.1 in [KMV] where no weight was considered. 
The proof of Theorem 2.1 is based on the following two lemmas, which are straightforward generalizations 
of Lemma 3.2 and Lemma 3.4 in [KMV] respectively.\par

\proclaim{Lemma 2.2} Let $u\in W_w^{1,q}(\Bbb R^n),$ $1<q\le 
n$ and $w(x)\in A_q(\Bbb R^n)$.  Then
$$\lim_{r\to 0} \intav_{B(x,r)} |u(y)-u(x)|^q\,w(y)\,dy=0\tag 2.3$$
except for $x$ in a set $E\subset \Bbb R^n$ of $(q,w)$-capacity zero. 
\endproclaim

In this paper we denote ${1\over {\int_{B(x,r)} w(y)\,dy}} \int_{B(x,r)}$ by $\intav_{B(x,r)}$.

The proof of this lemma is a straightforward adaptation of the proof of Lemma 3.2 in [MV1], replacing the 
initial definition in the proof of that lemma of $A_r u(x)$ by the new definition 
$$A_r u(x)=r^q\, \intav _{B(x,r)} |u(y)-u(x)|^q\,w(y)\,dy,$$
and using then Lemma 3.1 in [MV1] and the fact that smooth functions with compact support are dense 
in $W_w^{1,q}(\Bbb R^n)$ whenever the weight $w$ is in the class $A_q$ (see [K]).\par

\proclaim{Lemma 2.4} [HKM, Theorem 3.34] Let $u$ be a weighted $\Cal A$-harmonic  
function in $\Bbb B^n,$ and fix $1<q\le 
n$ and $w(x)\in A_q(\Bbb R^n)$. Then, there exists a constant $C$ such  
that for each ball $B=B(x,r)\subset \Bbb B^n$ and all $a\in \Bbb R$ 
$$\sup_{\frac 1 2 B}|u(y)-a|\le C\,\,(\intav_B |u(y)-a|^q\,w(y)\,dy)^{1/q},$$
where $\frac 1 2 B=B(x,r/2).$
\endproclaim


\demo{Proof of Theorem 2.1} 
Since $\int_{\Bbb B^n} |\nabla u(x)|^q\,w(x)\,dx<\infty$, it follows  
from the Poincar\'e inequality that $u\in W_w^{1,q}(\Bbb B^n)$. Hence, by  
standard extension theorems, we may assume that $u\in W_w^{1,q}(\Bbb R^n)$. 
We  show that $u$ has a nontangential limit for each $x\in \partial \Bbb B^n$ for  
which (2.3) holds. The claim then follows from Lemma 2.2. \par
\noindent Fix a $w\in\partial\Bbb B^n$ for which (2.3) holds.
We denote by $C(w)$ the Stolz cone at $w$  
with a fixed given aperture. Then we can find a constant $c_n\ge 1$, depending  
only on the aperture and $n$, such that for all $x\in C(w)$ 
$$|w-x|\le c_n (1-|x|).$$ \par
\noindent Pick $x\in C(w)$. Then, 
we have that $$B(x, (1-|x|)/2)\subset B(w,(c_n+\frac {1}{2})(1-|x|)),$$  
and hence, by Lemma 2.4,
$$
\align
|u(x)-u(w)|&\le C\,\,(\intav_{B(x, (1-|x|)/2)} 
|u(y)-u(w)|^q\,w(y)\,dy)^{1/q}\\
&\le C^\prime \,\,(\intav_{B(w,(c_n+\frac {1}{2})(1-|x|))}  
|u(y)-u(w)|^q\,w(y)\,dy)^{1/q}.
\endalign
$$
The claim follows by applying (2.3).
\qed
\enddemo


In Theorem 2.1, $\Cal A$-harmonicity was not essential but merely the version  
of the weak weighted Harnack inequality (Lemma 2.4) satisfied by the solutions to a  
large class of elliptic nonlinear P.D.E.'s.


\beginsection{ $\S$3.  Nontangential limits for Quasiregular Mappings.}

Let $W^{1,n}_{\text{loc}}(\Bbb B^n)$ denote the local Sobolev space of  
functions in $L^n_{\text{loc}}(\Bbb B^n)$ whose distributional derivatives  
belong to $L^n_{\text{loc}}(\Bbb B^n)$. Consider a mapping

$$f\colon \Bbb B^n\to \Bbb R^n$$
whose coordinate functions belong to $W^{1,n}_{\text{loc}}(\Bbb B^n)$. Denote  
by $J_f(x)$ the Jacobian determinant $\text {det}(Df(x))$. For a.e. $x\in \Bbb  
B^n$ the {\it dilatation} of $f$ is defined by
$$K(x)=\frac{|Df(x)|^n}{J_f(x)},$$
and it satisfies $K(x)\ge c_n$. If $K(x)\in L^{\infty}(\Bbb B^n),$  then $f$ is  
said to be a quasiregular mapping. \par
It is well known, see [HKM], that if $f$ is a nonconstant K-quasiregular  
mapping in $\Bbb B^n$ and $b\in \Bbb R^n$, the function $u(x)=\log|f(x)-b|$ is 
a weighted $\Cal A$-harmonic function in $\Bbb B^n\setminus f^{-1}(b)$ of type $q=n$, weight $w=1$, 
$\alpha={1\over K}$, and $\beta=K$. Therefore, all the results  
of the previous section apply to quasiregular mappings. \par
In this section we prove that a certain restriction on the  
growth of the multiplicity function of $f$ implies the existence of  
nontangential limits. \par
After these preliminaries, let us prove the following result.

\proclaim{Theorem 3.1} Let $f$ be a bounded quasiregular mapping of $\Bbb B^n$. Let $w$ be a nonnegative weight 
defined by
$$w(x)= \sum_{j=0}^\infty c_j\,|1-|x||^{q-1}\,\aleph_{R_j}(x),\tag 3.2$$
where $1<q<n$ and $R_j=\lbrace x\colon 1-2^{-j}<|x|\leq 1-2^{-j-1}\rbrace,$ $j=0,1,2,\ldots$, $\aleph_{R_j}$ is the 
characteristic function of $R_j$, and the $c_j$'s are positive constants satisfying that 
for some positive $b$ $\sum_{j=1}^\infty {{c_j}\over{j^{b\,{q\over n}}}}<\infty$. 
Let us assume that for the same positive $b$ we have that 
$$N(f, B(0,r))\leq C\,(1-r)^{-(n-1)}\,\Biggl( {1\over{\log 
\Bigl({1\over {1-r}}\Bigr)}}\Biggr)^{b}\tag 3.3$$ 
for all $0<r<1$.\par
Then 
$$\int_{\Bbb B^n(0,1)} |Df(x)|^q\,w(x)\,dw<\infty.$$
\endproclaim

It is important to observe that here we do not need to assume that the weight $w(x)$ as defined in (3.2) 
is in $A_q(\Bbb R^n)$.\par

\demo{Proof of Theorem 3.1} Let $R_j=\lbrace x\colon 1-2^{-j}<|x|\leq 1-2^{-j-1}\rbrace,$ $j=1,2,\ldots$. Now 

$$\sum_{j=1}^\infty \int_{R_j} |Df(x)|^q\,w(x)\,dx\leq C\sum_{j=1}^\infty (\int_{R_j}  
|Df(x)|^n\,dx)^{q/n}\, (\int_{R_j}  
w(x)^{{n\over{n-q}}}\,dx)^{{{n-q}\over n}},$$
and by a change of  variables [BI, 8.3] we arrive at

$$
\align
\sum_{j=1}^\infty \int_{R_j} |Df(x)|^q\,w(x)\,dx&\leq  C\sum_{j=1}^\infty  
(\int_{f(R_j)} n(y,f,R_j)\,dy)^{q/n}\\& (\int_{R_j}  
w(x)^{{n\over{n-q}}}\,dx)^{{{n-q}\over n}}.
\endalign
$$
Since

$$
\align
n(y,f,R_j)&\leq N(f,B(1-2^{-j-1}))\\
&\leq C\,(1-(1-2^{-j-1}))^{-(n-1)}\,\Biggl( {1\over{\log 
\Bigl({1\over {1-(1-2^{-j-1})}}\Bigr)}}\Biggr)^{b}\\
&=C\, 2^{j(n-1)}\,{1\over{j^b}},
\endalign
$$ 

\noindent and since $w(x)= \sum_{j=0}^\infty c_j\,|1-|x||^{q-1}\,\aleph_{R_j}(x)$, then in 
$R_j$, 
$$w(x)\leq C\, 2^{-j(q-1)}\,c_j.$$
Thus

$$\sum_{j=1}^\infty \int_{R_j} |Df(x)|^q\,w(x)\,dx\leq C\sum_{j=1}^\infty  
2^{j(n-1){q\over n}}\,{1\over{j^{b\,{q\over n}}}}\,2^{-j({{n-q}\over n})}\,
2^{-j(q-1)}\,c_j .$$
Thus, after simplifying we have that

$$
\align
\sum_{j=1}^\infty \int_{R_j} |Df(x)|^q\,w(x)\,dx& \leq C\sum_{j=1}^\infty {{c_j}\over{j^{b\,{q\over n}}}}
\\&=C\sum_{j=1}^\infty {{c_j}\over{j^{b\,{q\over n}}}}<\infty,
\endalign
$$
by assumption. Therefore $\int_{\Bbb B^n\setminus B(0,1/2)} |Df(x)|^q\,w(x)\,dx<\infty$, thus  
$$\int_{\Bbb B^n} |Df(x)|^q\,w(x)\,dx<\infty$$ 
and the claim will follow.
\qed
\enddemo
\vskip 0.2in

Let $w$ be a weight as in Theorem 3.1 with the constants $c_j$, 
$j=1,2,\ldots$ satisfying the hypothesis of the theorem. Let us assume at this stage that $w$ belongs to $A_q(\Bbb R^n)$. \par

It follows from [A, Theorem 6.1] that
$$(q,w)-\text{ cap} (\Bbb B^n(x,r))\approx \Biggl[\int_r^\infty t^{{{(q-n)}\over{(q-1)}}}\,
\intav_{\Bbb B(x,t)} w^{-{1\over {q-1}}}\,{{dt}\over t}\Biggr]^{1-q}.\tag 3.4$$
Since we are assuming that $w\in A_q$, we have that
$$\intav_{\Bbb B(x,t)} w^{-{1\over {q-1}}}\leq C\, {1\over{\Biggl(\intav_{\Bbb B(x,t)} w\Biggr)^{{1\over{q-1}}}}},$$
and thus dividing by something bigger, 

$$(q,w)-\text{ cap} (\Bbb B^n(x,r))\geq C\,\Biggl[\int_r^\infty t^{{{(q-n)}\over{(q-1)}}}\,
{1\over{\Biggl(\intav_{\Bbb B(x,t)} w\Biggr)^{{1\over{q-1}}}}}\,{{dt}\over t}\Biggr]^{1-q}.\tag 3.5$$
Using now the explicit formula for our weight $w$, we have that

$$w(\Bbb B(x,t))\approx \sum_{j=j_0}^\infty 2^{-j(n+q-1)}\, c_j
\approx t^{n+q-1}\,c_{j_0(t)},$$

\noindent where $j_0(t)$ is a positive integer satisfying that $t\approx 2^{-j_0(t)}$. By definition 

$$\intav_{\Bbb B(x,t)} w={1\over {t^n}}\,w(\Bbb B(x,t)),$$

\noindent then, substituting in (3.5) we have that

$$(q,w)-\text{ cap} (\Bbb B^n(x,r))\geq C\,\Biggl[\int_r^\infty t^{{{(q-n)}\over{(q-1)}}}\,
{1\over{\Biggl(t^{q-1}\,c_{j_0(t)}
\Biggr)^{{1\over{q-1}}}}}\,{{dt}\over t}\Biggr]^{1-q}.$$

Since for any $t$ between $r$ and $\infty$ we can assume without loss of generality that 

$${1\over{\Biggl(t^{q-1}\,c_{j_0(t)}
\Biggr)^{{1\over{q-1}}}}}\leq {1\over{\Biggl(r^{q-1}\, c_{j_0(r)}\Biggr)^{{1\over{q-1}}}}},$$
thus
$$(q,w)-\text{ cap} (\Bbb B^n(x,r))\geq C\,r^{q-1}\, c_{j_0(r)}\,
\Biggl[\int_r^\infty t^{{{(q-n)}\over{(q-1)}}}\,{{dt}\over t}\Biggr]^{1-q}.\tag 3.6$$
An easy computation shows that the right hand side of (3.6) is equal to 
$$C\,r^{q-1}\, c_{j_0(r)}\, r^{n-q}= C\,r^{n-1}\, c_{j_0(r)}.$$ 
Hence we have that,

$$(q,w)-\text{ cap} (\Bbb B^n(x,r))\geq C\,c_{j_0(r)}\,r^{n-1}.\tag 3.7$$

Let $E$ be a subset of $\partial \Bbb B^n(0,1)$ and denote the $(n-1)$ dimensional Hausdorff measure of $E$ by
$$\Lambda_{n-1}(E)=\lim_{\delta\to 0} \Bigl[ \inf \lbrace \sum_i r_i^{n-1}\colon E\subset \bigcup 
\Bbb B^n(x_i,r_i),\,0<r_i<\delta\rbrace\Bigr],$$
where the infimum is taken over all coverings of $E$ by balls of radii less than $\delta$. 

Let $\lbrace \Bbb B^n(x_i,r_i)\colon x_i\in E, \, 0\le r_i<\delta\rbrace$ be a covering of the set $E$. If we define by 
$$\Lambda_{n-1}^\delta(E)= \inf \lbrace \sum_i r_i^{n-1}\colon E\subset \bigcup 
\Bbb B^n(x_i,r_i),\,0<r_i<\delta\rbrace,$$ 
we have that $\Lambda_{n-1}(E)=
\lim_{\delta\to 0} \Lambda_{n-1}^\delta(E)$. Hence for any of those coverings,

$$\Lambda_{n-1}^\delta(E)\le \sum_i r_i^{n-1}< \sum_i r_i^{n-1}\,c_{j_0(r_i)},$$

\noindent provided the $r_i$ are small enough and since $c_{j_0(r)}$ goes to $\infty$ as $r$ goes to 0. 
It is clear from our construction that we can always choose our $x_i$'s such that (3.7) 
holds for all the balls $\Bbb B^n(x_i,r_i)$, and thus we have that
$$\Lambda_{n-1}^\delta(E)\le \,C\, 
\lbrace \sum_i (q,w)-\text{cap} (\Bbb B^n(x_i,r_i))\rbrace.\tag 3.8$$

Let $f$ be a bounded quasiregular mapping in $\Bbb B^n$ and $w$ be a weight as in Theorem 3.1. 
Let $E$ be the set of points in $\subset \partial \Bbb B^n(0,1)$ for which 
the nontangential limit of $f$ does not exist. If our weight $w$ will be in $A_q(\Bbb R^n)$, by Theorem 2.1 in section 
\S2 we will have that $(q,w)-\text{cap} (E)=0$.

Thus, by the definition of the weighted 
variational capacity we have that for any $\tilde\epsilon>0$ 
we can find a covering of $E$ such that $E\subset \bigcup_i \Bbb B^n(x_i,r_i)$ 
and 
$$\sum_i (q,w)-\text{cap} (\Bbb B^n(x_i,r_i)) \le (q,w)-\text{cap} (E)+\tilde\epsilon.\tag 3.9$$
Combining (3.8) and (3.9) we have that
$$\Lambda_{n-1}^\delta(E)\le C\,\Bigl[ (q,w)-\text{cap} (E)+\tilde\epsilon\Bigr]=C\,\tilde\epsilon$$
and since $\tilde \epsilon$ is arbitrary and independent of $\delta$, letting $\delta\to 0$ we will obtain that 
$\Lambda_{n-1}(E)=0$.
\vskip 0.2in

In our last two paragraphs we have made two assumptions in order to get our conclusion that 
$\Lambda_{n-1}(E)=0$ for the set $E\subset \partial \Bbb B^n$ where the nontangential limits of the mapping $f$ fail to 
exist. Namely, the multiplicity function of the mapping $f$ satisfies the following growth condition

$$N(f, B(0,r))\leq C\,(1-r)^{-(n-1)}\,\Biggl( {1\over{\log \Bigl({1\over {1-r}}\Bigr)}}\Biggr)^{b}$$
\noindent and our weight $w$ with the following explicit formula
 
$$w(x)= \sum_{j=0}^\infty c_j\,|1-|x||^{q-1}\,\aleph_{R_j}(x)$$
\noindent is in the Class $A_q$. Our next result shows that for our second assumption to be true, the constants 
$c_j$ have an exponential growth to infinity as $j\to \infty$. If that is the case, a simple 
computation that we will leave to the reader, will show that then we can not ascertain that
$$\int_{\Bbb B^n(0,1)} |Df(x)|^q\,w(x)\,dw$$
is finite, and then Theorem 2.1 can not be invoked. This argument shows that we have pushed our approach to the limit 
in the sense that Theorem 3.11 below is the best result we can obtain following our approach.

In [MV1], it was shown that the positive weight $w(x)=|1-|x||^\alpha$ for $x\in \Bbb R^n$ is in 
the class $A_q(\Bbb R^n)$ whenever $q>\alpha+1$. \par
The following argument shows that, in some sense, this can not be improved. Namely,\par

\proclaim{ Lemma 3.10} Let the positive weight $w$ be defined as in Theorem 3.1. If $w(x)$ belongs to the 
Muckenhoupt class $A_q(\Bbb R^n)$ after extending it to be symmetric outside the 
unit ball of $\Bbb R^n$, then the $c_j$'s are equivalent to 
$c_j=(1-|x|)^{-\epsilon}$ on each of the rings $R_j$, for some positive $\epsilon$, and hence our weight $w$ is 
equivalent to 
$$\sum_{j=0}^\infty |1-|x||^{q-1-\epsilon}\,\aleph_{R_j}(x),$$
for some positive $\epsilon$.
\endproclaim

This shows that if we require the weight $w(x)$ in Theorem 3.1 to be in the class $A_q$, this forces it to be 
equivalent to $|1-|x||^{q-1-\epsilon}$.\par

\demo {Proof of Lemma 3.10} Since the integrals that appear in the definition of the $A_q$-weights are invariant 
under rotations for the weights under consideration, it is enough to show that
$$\sup_{\Bbb B} \biggl( {1\over{|{\Bbb B}|}}\,\int_{\Bbb B} w(x)\,dx\biggr)\,\,\biggl( 
{1\over {|{\Bbb B}|}}\,\int_{\Bbb B} [w(x)]^{{1\over{1-q}}}\,
dx\biggr)^{q-1}<\infty, $$
where $\Bbb B$ is any ball in $\Bbb R^n$ whose center falls in the positive real axis. 
Moreover, we can assume that $x_0$, the center of the ball, is equal to $(1,0,\ldots,0)$.\par
Using polar coordinates, and letting $s=|x|$ we have 
 
$$
\align
A&=\biggl( {1\over{|{\Bbb B(x_0,r)}|}}\,\int_{\Bbb B(x_0,r)} w(x)\,dx\biggr)\,\,\biggl({1\over 
{|{\Bbb B(x_0,r)}|}}\,\int_{\Bbb B(x_0,r)} [w(x)]^{{1\over{1-q}}}\,dx\biggr)
^{q-1}\\
&={c\over{r^{nq}}}\,\biggl( \int_I\,\,\int_{\Bbb B(x_0,r)\cap S^{n-1}_s} w(s) 
\,dS\,ds\biggr)\biggl( \int_I\,
\int_{\Bbb B(x_0,r)\cap S^{n-1}_s} w(s)^{{1\over {1-q}}} \,dS\,ds\biggr)^{q-1},
\endalign
$$
where $I$ is an interval of length equivalent to $r$ on the positive real axis ending at 
$x_0$ and $c$ is a constant that depends only on $n$. 
Using the fact that 
$$\int_{\Bbb B(x_0,r)\cap S^{n-1}_s} dS\leq c\,\,r^{n-1}$$
we reduce the problem to one dimension, and thus we need to show that
$$A\leq {c\over{r^{q}}}\,\biggl( \int_I  w(s)\,ds\biggr)\biggl( \int_I w(s)^{{1\over {1-q}}} \,ds\biggr)^{q-1}<\infty.$$
It is clear from the definition of the weight $w$, that it is enough to show that 
the above integral is bounded by a constant independent of the length of the interval 
$I$ whenever $I$ is an interval ending at 1. That is, we need to show that

$$
\align
A&\leq {c\over{r^{q}}}\,\biggl( \int_{r_0}^1  \sum_{j=0}^\infty c_j\,(1-s)^{q-1}\,\aleph_{(1-2^{-j},1-2^{-j-1}]}(s)
\,ds\biggr)\\
&\biggl( \int_{r_0}^1 \Biggl(\sum_{j=0}^\infty c_j\,(1-s)^{q-1}\,\aleph_{(1-2^{-j},1­ 2{-j-1}]}(s)\Biggr)^{{1\over {1-q}}}
 \,ds\biggr)^{q-1}<\infty,
\endalign
$$
where $r=1-r_0$. Without loss of generality we can assume that $r_0=1-2^{-j_0}$. 
Letting $t=1-s$ and performing the linear change of variable we need to show that

$$
\align
A&\leq {c\over{r^{q}}}\,\biggl( \int_0^{1-r_0}  \sum_{j=j_0}^\infty c_j\,t^{q-1}\,
\aleph_{(2^{-j-1},2^{-j}]}(t)\,dt\biggr)\\
&\biggl( \int_0^{1-r_0} \Biggl(\sum_{j=j_0}^\infty c_j\,t^{q-1}\,
\aleph_{(2^{-j-1},2^{-j}]}(t)\Biggr)^{{1\over {1-q}}}
 \,dt\biggr)^{q-1}\\
 &\leq {c\over{r^{q}}}\,\biggl( \sum_{j=j_0}^\infty c_j\,\int_{2^{-j-1}}^{2^{-j}} t^{q-1}\,dt\biggr)
 \biggl( \sum_{j=j_0}^\infty {c_j}^{{1\over{1-q}}}\,\int_{2^{-j-1}}^{2^{-j}} {{dt}\over t}\,dt\biggr)^{q-1}\\
 &={c\over{r^{q}}}\,\biggl( \sum_{j=j_0}^\infty c_j\,2^{-j\,q}\, (1-2^{-q})\biggr)
 \biggl( \sum_{j=j_0}^\infty {c_j}^{{1\over{1-q}}}\,\biggr)^{q-1},
\endalign
$$

\noindent where the constants $c$ in the above chain of inequalities might be different from line to line, but in any 
case is independent of the $r$. We have also used the fact that $\int_{2^{-j-1}}^{2^{-j}} {{dt}\over t}\,dt=\ln 2$.\par

It is clear by looking at the two infinite series after the last equality, that in order for $A$ to be bounded 
by a constant independent of $r$ we need to choose $c_j=\eta^j$ for $1<\eta<2^q$. 
But, this choice of the $c_j$'s on each of the rings $R_j$ is equivalent to say that on those rings 
$c_j=(1-|x|)^{-\epsilon}$ for some positive $\epsilon$ as we wanted to show.
\qed
\enddemo

Lemma 3.10 shows that we can not construct a weight around the weight 
$|1-|x||^{q-1}$ of the form (3.2) which is still in $A_q$ and for which the sequence of $c_j$'s 
goes to infinity as $j\to\infty$ less than exponentially. \par

Finally, we will prove the main result in this paper that improves Theorem 4.1 in [KMV] and the results of 
Martio and Srebro in [MS], on the size 
of the sets on the boundary of the unit ball of $\Bbb R^n$ where the nontangential limits might not exist. 
We will state our result for bounded quasiregular mappings in $\Bbb B^n$, but the same conclusion will follow if 
we impose a growth condition on the mapping $f$ of the type $|f(x)|\leq C\, (1-|x|)^{-b}$ for some $0\leq b<\infty$, 
as the argument in the proof of Theorem 4.1 in [KMV] shows.\par
Let us state our main result in this paper.\par

\proclaim  {Theorem 3.11} Let $f$ be a bounded quasiregular mapping of $\Bbb B^n$, and suppose that, for some 
$0\leq a<n-1$,
$$N(f, B(0,r))\leq C\,(1-r)^{-a}$$ 
for all $0<r<1$. Then the set of points $x_0\in E\subset \partial \Bbb B^n(0,1)$ for which 
the nontangential limit of $f$ does not exist has Hausdorff dimension less 
than or equal to $a$, i.e. $\text{dim}_H(E)\leq a$.
\endproclaim

Before we pass to the proof of this Theorem, notice that for any $0< a<n-1$ we have that 
$a<{{na}\over{1+a}}$, which shows that our result improves the one in [KMV].

\demo{Proof of Theorem 3.11} Arguing as in the proof of Theorem 3.1 we arrive at

$$
\align
\sum_{j=1}^\infty \int_{R_j} |Df(x)|^q\,w(x)\,dx&\leq  C\sum_{j=1}^\infty  
(\int_{f(R_j)} n(y,f,R_j)\,dy)^{q/n}\\ 
& (\int_{R_j}  
w(x)^{{n\over{n-q}}}\,dx)^{{{n-q}\over n}}.
\endalign
$$
Since

$$
\align
n(y,f,R_j)&\leq N(f,B(1-2^{-j-1}))\\&\leq C\,(1-(1-2^{-j-1}))^{-a}\\
&=C\, 2^{j\,a}.
\endalign
$$ 

For any positive $\epsilon$ we choose $w(x)= |1-|x||^{\alpha}$ with $\alpha+1<q$. 
By Lemma 5.1 in [MV1] $w\in A_q(\Bbb R^n)$ and we have that in $R_j$, 
$$w(x)\leq C\, 2^{-j\,\alpha}.$$
Then we have that,

$$\sum_{j=1}^\infty \int_{R_j} |Df(x)|^q\,w(x)\,dx\leq C\sum_{j=1}^\infty  
2^{j\,a\,{q\over n}}\,\,2^{-j({{n-q}\over n})}\,
2^{-j\,\alpha}.$$

\noindent Thus, after simplifying we have that

$$
\align
\sum_{j=1}^\infty \int_{R_j} |Df(x)|^q\,w(x)\,dx&\leq C\sum_{j=1}^\infty 2^{j(\,a\,{q\over n}-({{n-q}\over n})
-\alpha)}\\&
=C\sum_{j=1}^\infty 2^{j((1+a)\,{q\over n}-(\alpha+1))}<\infty,
\endalign
$$

\noindent if and only if $(1+a)\,{q\over n}-(\alpha+1)$ is negative. That is, if and only if 
$(1+a)\,{q\over n}<\alpha+1<q$. Thus, for fixed $0\leq a<n-1$, let $\epsilon$ positive such that 
$q-1>\alpha=(1+a+\epsilon)\,{q\over n}-1>(1+a)\,{q\over n}-1$.

Therefore $\int_{\Bbb B^n\setminus B(0,1/2)} |Df(x)|^q\,w(x)\,dx<\infty$, thus  
$\int_{\Bbb B^n} |Df(x)|^q\,w(x)\,dx<\infty$ and the claim will follow.

Applying Theorem 2.1 to each of the components of the bounded quasiregular mapping $f$, 
we have that the mapping $f$ has nontangential limits on all radii terminating outside a set of $(q,w)$-capacity  
zero, where $1<q\leq n$ and the weight $w$ as above.\par

Next, we will show that the set $E$ has Hausdorff
 dimension less than or equal to $\alpha+n-p$. \par
 
We start by observing the following weak type estimate
$$\Lambda^{w,\infty}_s\Biggl(\lbrace y\in \Bbb R^n\colon M^w_{s,p} f(y)>t\rbrace \Biggr) \leq {{c(s,w)}
\over{t^p}}\,\int_{\Bbb R^n} |f(x)|^p\,w(x)\,dx,\tag 3.12$$
whose proof is actually contained in the proof of Lemma 3.1 in [MV1]. 
 Let $u$ be a function in $C^\infty_0(B(y,R))$, as in [HKM, Lemma 2.30] write 
$$|u(y)|\leq c\,\int_0^\infty {1\over{r^n}}\,\int_{B(y,r)} |\nabla u(x)|\,dx\,dr.$$
Inserting $w(x)^{1/p}\,\,w(x)^{-1/p}$ and applying H\"older's inequality and the $A_p$ condition for 
the weight $w$ one can easily adapt the proof of [HKM, Lemma 2.30] to obtain
$$|u(x)|\leq c\, R^{1+s/p}\, M_{s,p}^w [|\nabla u|(x)]\tag 3.13$$
for any $x\in \Bbb R^n$ as long as $s>-p$.\par
Combining (3.12) and (3.13) we obtain
$$\Lambda^{w,\infty}_s\Biggl(\lbrace x\in \Bbb B(y,R)\colon |u(y)|>t\rbrace \Biggr) \leq c\, 
{{R^{s+p}}\over {t^p}}\,\int_{\Bbb B(y,R)} |\nabla u(x)|^p\, w(x)\,dx.\tag 3.14$$
The proof of [HKM, Theorem 2.26] can be carried out in the weighted case by using the estimate 
(3.14), thus the $cap_{p,w}(E)=0$ implies that $\Lambda^{w,\infty}_s(E)=0$ for all $s>-p$. \par
It follows that given any positive $\epsilon$ one can cover the set $E\subset \bigcup_i \Bbb B(x_i,r_i)$, 
such that 
$$\sum_{i} r_i^{s}\,w(\Bbb B(x_i,r_i))<\epsilon.$$
Note that the centers of the balls $x_i$ may be taken on $\partial \Bbb B^n$. Using now the special nature 
of our weight we have that $w(\Bbb B(x_i,r_i))\approx r_i^{n+\alpha}$. Therefore, 
$$\sum_{i} r_i^{n+\alpha+s}<c\,\epsilon.$$
Using for example [HKM, Lemma 2.25] we conclude that the Hausdorff dimension of the set $E$ is less than 
or equal to $n+\alpha+s$ for any $s>-p$.
  for any $s>-p$.
Thus, $\text{dim}_H(E)\leq n+\alpha-q$.\par
  
  Finally, choosing $q=n$, since we have that 
  $q-1>\alpha=(1+a+\epsilon)\,{q\over n}-1>(1+a)\,{q\over n}-1$, then $\alpha=a+\epsilon$ and thus
  $\text{dim}_H(E)\leq a+\epsilon$. Letting $\epsilon$ go to zero, we obtain 
  the desired result and the Theorem is proved.
  \qed
  \enddemo

The main idea in the proof of Theorem 3.11 is to reduce our situation to the  
case of a weighted Dirichlet finite quasiregular mapping, and then apply the results of the section \S2.\par

\proclaim{Remarks 3.15} 1) Theorem 3.11 generalizes a result of Koskela, Manfredi and Villamor [KMV] 
stating that bounded quasiregular mapping $f$ in $\Bbb B^n$ satisfying the growth  
condition $N(f, B(0,r))\leq C\,(1-r)^{-a}$ for some $0\leq a<n-1$ 
on the multiplicity function has nontangential limits   
at all points on the boundary of the unit ball except possibly on a set whose 
Hausdorff dimension is strictly less than ${{n\,a}\over {1+a}}<n-1$, since for 
any $0<a<n-1$ we have that $a<{{n\,a}\over {1+a}}$. It also improves on a result of Martio and Srebro [MS], who prove 
Theorem 3.11 for locally injective bounded quasiregular mappings in the unit ball of $\Bbb R^n$. In that paper, Martio 
and Srebro construct explicit examples of locally injective bounded quasiregular mappings that 
show that both their and our results are sharp.\par
2) We will obtain the same conclusion as in Theorem 3.11 if rather than assuming that the mapping $f$ is bounded 
and with the same growth condition on its multiplicy function, we will assume that 
$|f(x)|\leq C\,(1-|x|)^{-b}$ for some $0<b<\infty$. Namely, the following theorem holds:
\endproclaim

\proclaim  {Theorem 3.16} Let $f$ be a bounded quasiregular mapping of $\Bbb B^n$, and suppose that, for some 
$0\leq a<n-1$,
$$N(f, B(0,r))\leq C\,(1-r)^{-a}$$ 
for all $0<r<1$. If $|f(x)|\leq C\,(1-|x|)^{-b}$ for some $0<b<\infty$, 
then the set of points $x_0\in E\subset \partial \Bbb B^n(0,1)$ for which 
the nontangential limit of $f$ does not exist has Hausdorff dimension less 
than or equal to $a$, i.e. $\text{dim}_H(E)\leq a$.
\vskip 0.2in
The proof of this theorem requires, as in the proof of Theorem 4.1 in [KMV], an initial 
modification of our mapping $f$ by composing it with the mapping $g(x)= x\,|x|^{\epsilon-1}$, and 
picking $\epsilon>0$ so that $h(x)=g(f(x))$ satisfies $|h(x)|\leq C\, (1-|x|)^{-s}$, for some $s\geq 0$ with 
$(1+a+n\,s)<n$.
\endproclaim


\bye
\end